%
\def\conv{\mathop{\vrule height2,6pt depth-2,3pt 
    width 5pt\kern-1pt\rightharpoonup}}

\advance\vsize by 1 true cm
%
\def\dess #1 by #2 (#3){
  \vbox to #2{
    \hrule width #1 height 0pt depth 0pt
    \vfill
    \special{picture #3} 
    }
  }

\def\dessin #1 by #2 (#3 scaled #4){{
  \dimen0=#1 \dimen1=#2
  \divide\dimen0 by 1000 \multiply\dimen0 by #4
  \divide\dimen1 by 1000 \multiply\dimen1 by #4
  \dess \dimen0 by \dimen1 (#3 scaled #4)}
  }
%
\def \trait (#1) (#2) (#3){\vrule width #1pt height #2pt depth #3pt}
\def \fin{\hfill
	\trait (0.1) (5) (0)
	\trait (5) (0.1) (0)
	\kern-5pt
	\trait (5) (5) (-4.9)
	\trait (0.1) (5) (0)
\medskip}
%


\font\sevenbf=cmbx7

\baselineskip=15pt
\abovedisplayskip=15pt plus 4pt minus 9pt
\belowdisplayskip=15pt plus 4pt minus 9pt
\abovedisplayshortskip=3pt plus 4pt
\belowdisplayshortskip=9pt plus 4pt minus 4pt
\let\epsilon=\varepsilon

\def\biblio #1 #2\par{\parindent=30pt\item{}\kern -30pt\rlap{[#1]}\kern
30pt #2\smallskip}
 %
\catcode`\@=11
\def\@lign{\tabskip=0pt\everycr={}}
\def\equations#1{\vcenter{\openup1\jot\displ@y\halign{\hfill\hbox
{$\@lign\displaystyle##$}\hfill\crcr
#1\crcr}}}
\catcode`\@=12
%
\def\pmb#1{\setbox0=\hbox{#1}%
\hbox{\kern-.04em\copy0\kern-\wd0
\kern.08em\copy0\kern-\wd0
\kern-.02em\copy0\kern-\wd0
\kern-.02em\copy0\kern-\wd0
\kern-.02em\box0\kern-\wd0
\kern.02em}}
%
\def\undertilde#1{\setbox0=\hbox{$#1$}
\setbox1=\hbox to \wd0{$\hss\mathchar"0365\hss$}\ht1=0pt\dp1=0pt
\lower\dp0\vbox{\copy0\nointerlineskip\hbox{\lower8pt\copy1}}}
%

%

\def\maj#1#2,{\rm #1\sevenrm #2\rm{}}
\def\Maj#1#2,{\bf #1\sevenbf #2\rm{}}
\outer\def\lemme#1#2 #3. #4\par{\medbreak
\noindent\maj{#1}{#2},\ #3.\enspace{\sl#4}\par
\ifdim\lastskip<\medskipamount\removelastskip\penalty55\medskip\fi}

\def\Remark #1. {\noindent{\Maj REMARK,\ \bf #1. }}

\outer\def\Lemme#1#2 #3. #4\par{\medbreak
\noindent\Maj{#1}{#2},\ \bf #3.\rm\enspace{\sl#4}\par
\ifdim\lastskip<\medskipamount\removelastskip\penalty55\medskip\fi}



\def\Notation #1. {\noindent{\Maj NOTATION,\ \bf #1. }}

\def\Example #1. {\noindent{\Maj EXAMPLE,\ \bf #1. }}

\hfuzz=1cm


\catcode`\ˆ=\active     \def ˆ{\`a}
\catcode`\‰=\active     \def ‰{\^a}
\catcode`\=\active     \def {\c c}
\catcode`\Ž=\active    \def Ž{\'e} 

\catcode`\=\active   \def {\`e}
\catcode`\=\active   \def {\^e}
\catcode`\'=\active   \def '{\"e}
\catcode`\"=\active   \def "{\^\i}
\catcode`\•=\active   \def •{\"\i}
\catcode`\™=\active   \def ™{\^o}
\catcode`\š=\active   \defš{}
\catcode`\=\active   \def {\`u}
\catcode`\ž=\active   \def ž{\^u}
\catcode`\Ÿ=\active   \def Ÿ{\"u}
\catcode`\ =\active   \def  {\tau}
\catcode`\¡=\active   \def ¡{\circ}
\catcode`\¢=\active   \def ¢{\Gamma}
\catcode`\¤=\active   \def ¤{\S\kern 2pt}
\catcode`\¥=\active   \def ¥{\puce}
\catcode`\§=\active   \def §{\beta}
\catcode`\¨=\active   \def ¨{\rho}
\catcode`\©=\active   \def ©{\gamma}
\catcode`\­=\active   \def ­{\neq}
\catcode`\°=\active   \def °{\ifmmode\ldots\else\dots\fi}
\catcode`\±=\active   \def ±{\pm}
\catcode`\²=\active   \def ²{\le}
\catcode`\³=\active   \def ³{\ge}
\catcode`\µ=\active   \def µ{\mu}
\catcode`\¶=\active   \def ¶{\delta}
\catcode`\·=\active   \def ·{\Sigma}
\catcode`\¸=\active   \def ¸{\Pi}
\catcode`\¹=\active   \def ¹{\pi}
\catcode`\»=\active   \def »{\Upsilon}
\catcode`\¾=\active   \def ¾{\alpha}
\catcode`\À=\active   \def À{\cdots}
\catcode`\Â=\active   \def Â{\lambda}
\catcode`\Ã=\active   \def Ã{\sqrt}
\catcode`\Ä=\active   \def Ä{\varphi}
\catcode`\Å=\active   \def Å{\xi}
\catcode`\Æ=\active   \def Æ{\Delta}
\catcode`\Ç=\active   \def Ç{\cup}
\catcode`\È=\active   \def È{\cap}
\catcode`\Ï=\active   \def Ï{\oe}
\catcode`\Ñ=\active   \def Ñ{\to}
\catcode`\Ò=\active   \def Ò{\in}
\catcode`\Ô=\active   \def Ô{\subset}
\catcode`\Õ=\active   \def Õ{\superset}
\catcode`\Ö=\active   \def Ö{\over}
\catcode`\×=\active   \def ×{\nu}
\catcode`\Ù=\active   \def Ù{\Psi}
\catcode`\Ú=\active   \def Ú{\Xi}
\catcode`\Ü=\active   \def Ü{\omega}
\catcode`\Ý=\active   \def Ý{\Omega}
\catcode`\ß=\active   \def ß{\equiv}
\catcode`\à=\active   \def à{\chi}
\catcode`\á=\active   \def á{\Phi}
\catcode`\ä=\active   \def ä{\infty}
\catcode`\å=\active   \def å{\zeta}
\catcode`\æ=\active   \def æ{\varepsilon}
\catcode`\è=\active   \def è{\Lambda}  
\catcode`\é=\active   \def é{\kappa}
\catcode`\ë=\active   \defë{\Theta}
\catcode`\ì=\active   \defì{\eta}
\catcode`\í=\active   \defí{\theta}
\catcode`\î=\active   \defî{\times}
\catcode`\ñ=\active   \defñ{\sigma}
\catcode`\ò=\active   \defò{\psi}

\def\date{\number\day\
\ifcase\month \or janvier \or f\'evrier \or mars \or avril \or mai \or juin \or juillet \or ao\^ut  \or
septembre \or octobre \or novembre \or d\'ecembre \fi
\ \number\year}

\font \Ggras=cmb10 at 12pt
\font \ggras=cmb10 at 11pt

\def\sym{\fam\comfam\com}
\font\tensym=msbm10
\font\sevensym=msbm7
\font\fivesym=msbm5
\newfam\symfam
\textfont\symfam=\tensym
\scriptfont\symfam=\sevensym
\scriptfont\symfam=\fivesym
\def\sym{\fam\symfam\relax}

\def\R{{\sym R}}


\font \Ggras=cmb10 at 14pt
\font \ggras=cmb10 at 16pt

\def\biblio #1 #2\par{\parindent=30pt\item{}\kern -30pt\rlap{[#1]}\kern 30pt #2\smallskip}
\def\biblio #1 #2\par{\parindent=30pt\item{[]}\kern -30pt\rlap{[#1]}\kern 30pt #2\smallskip}

\bigskip
\centerline{\Ggras  ASYMPTOTIC BEHAVIOR OF STRUCTURES}
\centerline{\Ggras  MADE OF CURVED RODS} 
\vskip 5mm
\centerline{\tenrm {G. G{\sevenrm RISO} }}
\vskip 5mm
\centerline{Laboratoire J.-L. Lions--CNRS, Bo\^\i te courrier 187, Universit\'e  Pierre et
Marie Curie,}
\centerline{4~place Jussieu, 75005 Paris, France, Email: griso@ann.jussieu.fr}
\vskip 5mm

\noindent{\ggras \bf Abstract. } {\sevenrm  In this paper  we study  the asymptotic behavior of a structure made of
curved rods of thickness $\scriptstyle 2\delta$   when $\scriptstyle\delta\to 0$. This study is carried on within the frame of linear
elasticity by using the unfolding method. It is based on several decompositions of the structure displacements and on the passing
to the limit in fixed domains.

We show that any displacement of a structure is the sum of an elementary rods-structure displacement  (e.r.s.d.)
 concerning the   rods cross sections and  a residual one   related  to the deformation of the cross-section. The e.r.s.d. 
coincide with rigid body displacements in the junctions. Any e.r.s.d. is given by two functions belonging to  $\scriptstyle 
H^1({\cal S};\R^3)$ where $\scriptstyle {\cal S}$ is the skeleton structure (i.e. the set of the rods middle lines). One
of this function $\scriptstyle {\cal U}$ is the skeleton displacement, the other $\scriptstyle {\cal R}$ gives the cross-sections
rotation. We show that $\scriptstyle {\cal U}$ is the sum of an extensional displacement and an inextensional one. We establish a
priori estimates and then we characterize the unfolded limits of the rods-structure  displacements.

Eventually we pass to the limit  in the linearized elasticity system and using all results in [5],  on the one hand we obtain   a
variational  problem that is satisfied by the  limit extensional displacement,  and on the other hand,  a   variational  problem 
coupling the limit of   inextensional  displacement and the limit of the rods torsion angles.  } 
\smallskip
\noindent{\ggras \bf R\'esum\'e. } {\sevenrm  On \'etudie dans cet article le comportement asymptotique d'une structure
form\'ee de poutres courbes d'\'epaisseur $\scriptstyle 2\delta$  lorsque $\scriptstyle\delta\to 0$. Cette \'etude est men\'ee
dans le cadre de l'\'elasticit\'e lin\'eaire en utilisant la m\'ethode de l'\'eclatement.   Elle est bas\'ee sur plusieurs d\'ecompositions
des d\'eplacements de la structure, et  sur le passage \`a la limite dans des domaines fixes.  

On montre que tout d\'eplacement de la structure est la somme d'un d\'eplacement  \'el\'ementaire de structure-poutres
(d.e.s.p.) concernant  les sections droites des poutres  et d'un  d\'eplacement  r\'esiduel  li\'e aux  d\'eformations  de ces
sections droites. Les d.e.s.p.   co\"\i  ncident avec des d\'eplacements rigides dans les jonctions.  Tout
d.e.s.p.  est donn\'e par deux fonctions  appartenant ˆ $\scriptstyle H^1({\cal S};\R^3)$ o\`u $\scriptstyle {\cal S}$ est le
squelette de la structure (l'ensemble des lignes moyennes  des poutres). L'une de ces fonctions $\scriptstyle {\cal U}$ est le
d\'eplacement du squelette, l'autre $\scriptstyle {\cal R}$ donne la rotation des sections droites.  On montre  que $\scriptstyle
{\cal U}$   est la  somme  d'un d\'eplacement  extensionnel  et d'un d\'eplacement   inextensionnel.  On \'etablit des estimations
a  priori puis on caract\'erise  les limites des \'eclat\'es des  d\'eplacements des poutres de la structure.

 Pour finir on passe  \`a la  limite  dans le syst\`eme de l'\'elasticit\'e lin\'eaire, on obtient en utilisant tous les rŽsultats de
[5] d'une part  le   probl\`eme variationnel  v\'erifi\'e par la limite  des  d\'eplacements  extensionnels,   et d'autre part  le  
probl\`eme variationnel  couplant  la limite  des  d\'eplacements    inextensionnels et la limite des angles de torsion des poutres.}

\vfill\eject
\noindent {\ggras  1. Introduction} 
\medskip
This paper follows the one entitled ``Asymptotic behavior of curved rods by the unfolding method''.  Here we are interested in
the asymptotic behavior of  a  structure made of  curved rods in the framework of the linear elasticity according to the unfolding
method.  It consists in obtaining some displacements decompositions and then in passing to the limit in fixed domains.

The structure is made of curved rods  whose cross sections are discs of  radius $\delta$. The middle lines of the rods are
regular arcs (the skeleton  structure ${\cal S}$). Let us take a displacement of the  structure. According to the results obtained in
[5] we write its restriction  to each rod as the sum of a rod elementary displacement and a residual displacement.
But this family of rod elementary displacements  is not necessarily the restriction of a   $H^1$ displacement of the
structure because such displacements do not inevitably coincide in the junctions. This is the reason why we had to change the rod
elementary displacements and to replace them by rigid body displacements in the junctions. Hence we obtain an elementary
rods-structure displacement (e.r.s.d.) which has become an admissible displacement of the whole domain. Then, any displacement
of the structure is the sum of an e.r.s.d. and a residual one. An e.r.s.d. is linear in the cross section and is given by two functions
belonging to $H^1({\cal S};\R^3)$. The first one ${\cal U}$ is the skeleton displacement  and the second one gives
the rotation of the cross section. 
\vskip 1mm
 In order to account for the asymptotic behavior of the deformations tensor and the strains tensor we shall decompose the first
component ${\cal U}$ of the elementary  rods-structure displacement  into the sum of an extensional displacement and an
inextensional one.  An extensional displacement modifies the length of the middle lines while an inextensional displacement
does not change this length in a first approximation.
\vskip 1mm
Very few articles have been dedicated to the study of the junction of two elastic bodies. The case of the junction of a
three dimensional domain and a two dimensional one is explored in [1]  by P.G. Ciarlet et al. The junction of two plates is studied
in [4] and [9], and [11] deals with the junction of beams and plates. Concerning the junction of rods, we gave  in [5] a thorough
study of the structures made of straight rods. The asymptotic behavior is under the form of a coupling of the inextensional
displacements and the rods torsion angles. Le Dret gave a first study of the junction of two straight rods [8].  He starts his study
from the three-dimensional problem of linearized elasticity and uses a standard thin domain technique  (the rods are transformed in
fixed domain). Le Dret obtains the variational problem coupling the flexion displacements in both rods and some junction
conditions. 
\vskip 1mm
This paper is organized as follows : in Section 2 we desciribe the skeleton and the structure. In Section 3 we recall the definition
of an elementary rod displacement and we give the definition of an elementary rods-structure displacement (e.r.s.d.).
Lemma 3.2 is technical result which allows us to associate an e.r.s.d. $U_e$ to any displacement $u$ of the structure.  Then,
Lemme 3.4 gives this  e.r.s.d. estimates and those of the  difference $u-U_e$ for appropriate norms. This estimates have an
essential  importance in ouer study : they replace the  Korn inequality. In Section 4 we introduce the inextensional
displacements and  the extensional displacements. All these decompositions are introduced tofacilitate the study of the
asymptotic behaviors of displacements sequences (Section 5) and of the strain and stress tensors sequences (see [5]). 
In Section 5 we also give the set of the limit inextensional displacements.  In Section 6 we pose the problem of elasticity.
Thanks to theorems in  [5] we deduce  Theorem 6.2 which on the one hand gives us the variational problem verified by the
extensional displacements limit  and, on the other hand, the variational problem coupling the inextensional displacements 
limit and the rods torsion angles limit.
\vskip 1mm
 In this work we use the same notation as in [5]. The constants appearing in the estimates will always be independent from
$\delta$.  As a rule the Latin index $i$ takes values in $\{1,\ldots,N\}$ and the Latin indices $h$, $j$, $k$ and  $l$  take values in
$\{1, 2, 3\}$. We also use the Einstein convention of summation over repeated indices.
\bigskip
\noindent {\ggras 2. The structure made of curved rods} 
\medskip
 The Euclidian space $\R^3$ is related to the frame  $(O;\vec e_1,\vec e_2,\vec e_3)$. We denote by $x$ the running point  of
$\R^3$,  by $\|\cdot\|_{\R^3}$ the euclidian norm   and by $\cdot$ the scalar product in $\R^3$.  

For any open set $\omega$ in $\R^3$, and any displacement $u$ belonging to $H^1(\omega;\R^3)$, we put
$${\cal E}(u,\omega)=\int_{\omega}\gamma_{kl}(u)\gamma_{kl}(u),\quad \gamma_{kl}(u)={1\over 2}\Bigl\{{\partial u_k\over
\partial x_l}+{\partial u_l\over \partial x_k}\Bigr\}, \quad {\cal D}(u,\omega)=\int_{\omega}{\partial u_k\over\partial
x_l}{\partial u_k\over \partial x_l}.$$

Let there be given a family $\gamma_i$, $i\in\{1,\ldots,N\}$  of arcs. Each curve is parametrized by its arc length $s_i$,
$\overrightarrow{OM}(s_i)=\phi_i(s_i)$,
$i\in\{1,\ldots,N\}$. The mapping $\phi_i$ belongs to  ${\cal C}^3([0,L_i]\,;\R^3)$ and the arc
$\gamma_i$ is the range of $\phi_i$. The restriction of $\phi_i$ to interval $]0,L_i[$ is imbedded, the arcs can be closed.

\noindent The structure skeleton ${\cal S}$ is the curves set $\gamma_i$. 
\medskip
\Lemme HYPOTHESES 2.1.  We assume the following hypotheses on ${\cal S}$:

$\bullet$ ${\cal S}$ is connected,

$\bullet$ two arcs of ${\cal S}$ have their intersection reduced to a finite number of points, common points to two arcs are called
knots and the set of knots is denoted ${\cal N}$,

$\bullet$ the arcs are not tangent in one knot,

$\bullet$  Frenet frames
$(M_i(s);\overrightarrow{T}_i(s),\overrightarrow{N}_i(s),\overrightarrow{B}_i(s))$ are defined at each point of $[0,L_i]$. 

\noindent For any arc $\gamma_i$ the curvature $c_i$ is given by  Frenet formulae (see [5])
$$\left\{\eqalign{  
{d\overrightarrow{OM}_i\over ds_i}&=\overrightarrow{T}_i\qquad\quad ||\overrightarrow{T}_i||_{\R^3}=1\cr
{d\overrightarrow{T}_i\over ds_i}&=c_i\overrightarrow{N}_i\qquad||\overrightarrow{N}_i||_{\R^3}=1\qquad
\overrightarrow{B}_i=\overrightarrow{T}_i\land\overrightarrow{N}_i\cr }\right.\leqno(1)$$    Let us now introduce  the mapping 
$\Phi_i: [0,L_i]\times
\R^2\longrightarrow \R^3$ defined by 
$$\Phi_i(s_i,y_2,y_3)=\overrightarrow{OM}_i(s_i)+y_2\overrightarrow{N}_i(s_i)+y_3\overrightarrow{B}_i(s_i)\leqno(2)$$
\noindent There exists $\delta_0>0$ depending only on ${\cal S}$, such that  the restriction of $\Phi_i$ to the compact set
$[0,L_i]\times \overline{D}(O;\delta_0)$ is a ${\cal C}^1$-diffeomorphism of that set onto its range (see [5]).
 \medskip
\Lemme DEFINITION 2.2. The curved rod ${\cal P}_{\delta, i}$ with  rod center line $\gamma_i$ is defined as follows:
$${\cal P}_{\delta, i}=\Phi_i\bigl(]0,L_i[\times D(O; \delta)\bigr),\qquad \hbox{for}\; \delta\in]0,\delta_0].\leqno(3)$$ 
\noindent The whole structure ${\cal S}_\delta$ is
$${\cal S}_\delta=\bigcup_{i=1}^N{\cal P}_{\delta, i}\leqno(4)$$

  The   running  point of the cylinder $\omega_{\delta,i}=]0,L_i[\times D(O;\delta)$ is denoted $(s_i,y_2,y_3)$. The reference
domain linked to ${\cal P}_{\delta, i}$ is the open set $\omega_i=]0,L_i[\times D(O;1)$ obtained by transforming
$\omega_{\delta,i}=]0,L_i[\times D(O;\delta)$ by orthogonal  affinity of ratio $\displaystyle {1/ \delta}$. The 
 running  point of the reference domain $\omega_i$ is denoted $(s_i,Y_2,Y_3)$.
\vfill\eject
\noindent{\ggras 3. Elementary rods-structure displacements.}

\noindent The $H^1$ class displacements of ${\cal S}$ make up a space denoted
$$H^1({\cal S};\R^3)=\Bigl\{V\in \prod_{i=1}^N H^1(0,L_i;\R^3)\;| \; \forall A\in {\cal N},\;
 A\in\gamma_i\cap\gamma_j,\; V_i(a_i)=V_j(a_j)\Bigr\}\leqno(5)$$ where $a_i$   is the arc length of $A$ onto $\gamma_i$ 
and where $V_i=V_{|\gamma_i}$, $V_i(a_i)$ is denoted $V(A)$.

\noindent In the neighborhood of a knot  two rods or more can join together. We call junction the union of the parts  that are
common to two rods  at least.

\noindent There exists a real $\rho$ greater than or equal to 1, which depends only on
${\cal S}$, such that for any knot $A$. The rods junction at  $A$ be contained in the domain
$$A^\rho_\delta=\bigcup_{i=1}^N\Phi_i\bigl(]a_i-\rho\delta,a_i+\rho\delta[\times D(O;\delta)\bigr)$$ 

  We recall that an elementary displacement  of the  rod ${\cal P}_{\delta,i}$ (see [5]) is an element $\eta_i$ of
$H^1(\omega_{\delta,i};\R^3)$ that is written in the form
$$\eta_i(s_i,y_2,y_3)={\cal A}(s_i)+{\cal B}(s_i)\land\bigl(y_2\overrightarrow{N}_i(s_i) +y_3
\overrightarrow{B}_i(s_i)\bigr),\quad (s_i,y_2,y_3)\in
\omega_{\delta, i},$$ where ${\cal A}$ and ${\cal B}$ are elements of $H^1(0,L_i;\R^3)$. 
 
Let $u$ be a displacement belonging to $H^1({\cal S}_\delta;\R^3)$. We have shown (Theorem 4.3. in [5]) that there exist
elementary rod displacements $U^{'}_{e,i}$ such that
$$\left\{\eqalign{  &U^{'}_{e,i}(s_i,y_2,y_3)={\cal U}^{'}_i(s_i)+{\cal R}^{'}_i(s_i)\land\bigl(y_2\overrightarrow{N}_i(s_i) +y_3
\overrightarrow{B}_i(s_i)\bigr),\quad (s_i,y_2,y_3)\in
\omega_{\delta, i},\cr    &\vphantom{{C\over \delta^2}}{\cal D}(u-U^{'}_{e,i},\omega_{\delta,i})\le C\;{\cal E}(u,{\cal
P}_{\delta,i}),\qquad\|u-U^{'}_{e,i}\|^2_{ L^2(\omega_{\delta,i};\R^3)}\le C\;\delta^2{\cal E}(u,{\cal P}_{\delta,i}),\cr
&\delta^2\Bigl\|{d{\cal R}^{'}_i\over ds_i}\Bigr\|^2_{L^2(0,L_i;\R^3)}+\Bigl\|{d{\cal U}^{'}_i\over ds_i}  -{\cal
R}_i\land\overrightarrow{T}_i\Bigr\|^2_{L^2(0,L_i;\R^3)}\le {C\over
\delta^2}\;{\cal E}(u,{\cal P}_{\delta,i}).\cr}\right.\leqno(6)$$ The constants  depend only  on
the mid-lines $\gamma_i$, $i\in\{1,\ldots,N\}$. 

\Lemme DEFINITION 3.1. An {\bf elementary rods-structure displacement} is a displacement belonging to $H^1({\cal
S}_\delta;\R^3)$  whose restriction 
 to each rod is an elementary displacement and whose restriction to each junction is a rigid body displacement.

\noindent An elementary  rods-structure displacement depends of two functions belonging to $H^1({\cal S};\R^3)$.
\vskip 1mm
\noindent  In Lemma 3.3 we show that any displacement of  $H^1({\cal S}_\delta;\R^3)$  can be approximated by an elementary
rods-structure displacement. The Lemma 3.2 will allow us to construct such  displacement. In this lemma we consider a curved rod
${\cal P}_{\delta,i}$ of the structure ${\cal S}_\delta$, wich we denote without index $i$ to simplify.
\vskip 1mm
\Lemme LEMMA 3.2.   Let $u$ be a displacement  of  $H^1({\cal P}_\delta;\R^3)$, $A$ be a point of the center line of arc length
$a$, and let $\rho$  be a real greater or equal to 1,  fixed. There exists an elementary  rod displacement $U_e$, rigid in the domain
${\cal P}_{a,\delta}= \Phi\bigl((]a-\rho\delta,a+\rho\delta[\times D(O;\delta)\bigr)$, verifying
$$\left\{\eqalign{ & U_e(s, y_2,y_3)={\cal U}(s)+{\cal R}(s)\land\bigl(y_2\overrightarrow{N}(s) +y_3
\overrightarrow{B}(s)\bigr),\enskip(s,y_2,y_3)\in]0,L[\times D(O;\delta),\cr & {\cal D}(u-U_e,{\cal P}_\delta)\le C{\cal
E}(u,{\cal P}_\delta)\qquad ||u-U_e||^2_{L^2({\cal P}_\delta;\R^3)}\le C\delta^2{\cal E}(u,{\cal P}_\delta)\cr
&\delta^2\Bigl\|{d{\cal R}\over ds}\Bigr\|^2_{L^2(0,L;\R^3)}+\Bigl\|{d{\cal U}\over ds}  -{\cal
R}\land\overrightarrow{T}\Bigr\|^2_{L^2(0,L;\R^3)}\le {C\over \delta^2} {\cal E}(u,{\cal P}_\delta)\cr}\right.\leqno(7)$$

\noindent{\bf Proof : } Let ${\cal P}^{*}_{a,\delta}$ the curved rod portion
$\Phi\bigl((]a-(\rho+1)\delta,a+(\rho+1)\delta[\times D(O;\delta)\bigr)$ and $r$ be the rigid body displacement defined by
using the means of  $u(x)$ and  $x\land u(x)$ in the ball
$\displaystyle B(A;{\delta/ 5})$, (see Lemma 3.1 of [5]). The domain ${\cal P}^{*}_{a,\delta}$ is the union of a finite number
(which depends on $\rho$ and the center line) of bounded open sets, starshaped with respect to a ball. From various
implementations of Lemma 3.1 of [5] we derive the following estimates:
$${\cal D}(u-r;{\cal P}^{*}_{a,\delta})\le C{\cal E}(u;{\cal P}^{*}_{a,\delta})\qquad ||u-r||^2_{L^2({\cal
P}^{*}_{a,\delta};\R^3)}\le C{\cal E}(u;{\cal P}^{*}_{a,\delta})\leqno(8)$$ where constant $C$ depends  only on the center line
and $\rho$.

We modify now the elementary displacement given by $(6)$. The new elementary displacement $U_e$ is defined in the rod portion
${\cal P}_{a,\delta}$, so as to coincide, within this domain, with the rigid body displacement
$r$. We put ($M=\Phi(s,y_2,y_3)\in{\cal P}_{a,\delta}$) 
$$\left\{\eqalign{ U_e(s,y_2,y_3)&=r(M)=\vec a+\vec b\land\overrightarrow{AM}\cr &=\vec a+\vec
b\land\overrightarrow{M(a)M(s)}+\vec b\land\bigl(y_2\overrightarrow{N}(s)+y_3\overrightarrow{B}(s)\bigr)\cr (s,y_2,y_3)\in
&]a-\rho\delta,a+\rho\delta[\times D(O;\delta)\cr}\right.\leqno(9)$$  
\noindent  Let $m$ be  an even function $m$ belonging to ${\cal C}^\infty(\R;[0,1])$, which satisfies
$$m(t)=0\qquad\forall t\in [0,\rho],\qquad m(t)=1\qquad\forall t\in [\rho+1,+\infty[\leqno(10)$$ Let us take up the rod
elementary displacement defined by $(6)$. The new elementary displacement $U_e$ is defined by
$$\left\{\eqalign{   &{\cal U}(s)=m\bigl({s-a\over\delta}\bigr){\cal U}^{'}(s)+\Bigl(1-m\bigl({s-a\over\delta}
\bigr)\Bigr)\Bigl\{\vec a+\vec b\land\overrightarrow{M(a)M(s)}\Bigr\}\cr  &{\cal R}(s)=m\bigl({s-a\over\delta}\bigr){\cal
R}^{'}(s)+\Bigl(1-m\bigl({s-a\over\delta} \bigr)\Bigr)\vec b\cr &U_e(s,y_2,y_3)={\cal U}(s)+{\cal
R}(s)\land\bigl(y_2\overrightarrow{N}(s)+  y_3\overrightarrow{B}(s)\bigr)\quad (s,y_2,y_3)\in
\omega_\delta\cr}\right.\leqno(11)$$
\noindent Estimates $(6)$ and $(8)$ give us $(7)$.\fin

\noindent{\bf Remark 3.3. : } The former rod elementary displacement $U^{'}_e$ and the new one
$U_e$ satisfy the estimates
$$||{\cal U}-{\cal U}^{'}||^2_{L^2(0,L;\R^3)}+\delta^2\Bigl\|{d{\cal U}\over ds}-{d{\cal U}^{'}\over ds}\Bigr\|^2_{
L^2(0,L;\R^3)}+\delta^2||{\cal R}-{\cal R}^{'}||^2_{L^2(0,L;\R^3)}\le C{\cal E}(u,{\cal P}^{*}_{a,\delta})\leqno(12)$$ The
constant depends only on the rod center line and on $\rho$.\fin
\Lemme LEMMA 3.4.  Let $u$ be a displacement of $H^1({\cal S}_\delta;\R^3)$. There exists an elementary rods-structure 
displacement $U_e$ which is written in $\omega_{\delta,i}$ 
$$U_{e,i}(s_i,y_2,y_3)={\cal U}_i(s_i)+{\cal R}_i(s_i)\land\bigl(y_2\overrightarrow{N}_i(s_i)
+y_3\overrightarrow{B}_i(s_i)\bigr)$$ where
${\cal U}$ and ${\cal R}$ belongs to $H^1({\cal S};\R^3)$, such that
$$\left\{\eqalign{ &{\cal D}(u-U_e,{\cal S}_\delta)\le C{\cal E}(u,{\cal S}_\delta)\qquad ||u-U_e||^2_{L^2({\cal
S}_\delta;\R^3)}\le C\delta^2{\cal E}(u,{\cal S}_\delta)\cr &\sum_{i=1}^N\delta^2\Bigl\|{d{\cal R}_i\over
ds_i}\Bigr\|^2_{L^2(0,L_i;\R^3)} +\sum_{i=1}^N \Bigl\|{d{\cal U}_i\over ds_i} -{\cal R}_i\land\overrightarrow{T}_i
\Bigr\|^2_{L^2(0,L_i;\R^3)}\le {C\over \delta^2}{\cal E}(u,{\cal S}_\delta)\cr}\right.\leqno(13)$$

\noindent{\bf Proof : } Let us take a knot $A$. The ball centered in $A$ and of radius
 $\displaystyle {\delta / 5}$ is included into the junction $A^\rho_\delta$. Lemma 3.2 gives us then an elementary displacement
of each rod joining in $A$, rigid in $A^\rho_\delta$.  The estimates of Lemma 3.4. are immediate consequences of $(7)$. \fin  
\noindent The skeleton ${\cal S}$ is clamped at some arcs ends. Let  $\Gamma_0$ be the set of
${\cal S}$ clamping points. The whole structure ${\cal S}_\delta$ is then clamped onto domains that are discs. We denote
$\Gamma_{0\delta}$ the set of clamped points of ${\cal S}_\delta$. The space of admissible displacements of ${\cal S}$
(respectively ${\cal S}_\delta$) is denoted  $H^1_{\Gamma_0}({\cal S};\R^3)$ (respectively  $H^1_{\Gamma_0}({\cal
S}_\delta;\R^3)$).
$$\left\{\eqalign{ H^1_{\Gamma_0}({\cal S};\R^3)&=\Bigl\{u\in H^1({\cal S};\R^3)\ \ | \ \
u=0\quad\hbox{on}\enskip\Gamma_0\Bigr\}\cr H^1_{\Gamma_0}({\cal S}_\delta;\R^3)&=\Bigl\{u\in H^1({\cal
S}_\delta;\R^3)\ \ | \ \ u=0\quad\hbox{on}\enskip\Gamma_{0\delta}\Bigr\}\cr}\right.\leqno(13)$$
\noindent If a displacement $u$ belongs to $H^1_{\Gamma_0}({\cal S}_\delta;\R^3)$, the elementary rods-structure
displacement $U_e$ given by Lemma 3.4 is also an admissible displacement which belongs to $H^1_{\Gamma_0}({\cal
S}_\delta;\R^3)$.
\medskip
\noindent{\ggras 4. The  inextensional and extensional displacements } 
\medskip
\noindent The space of admissible displacements of ${\cal S}$ is equipped with the inner product 
$$<U,V>=\displaystyle\sum_{i=1}^N\int_0^{L_i}{dU_i\over ds_i}\cdot{dV_i\over ds_i}\leqno(14)$$ This inner product turns
$H^1_{\Gamma_0}({\cal S};\R^3)$ into a Hilbert space.

\Lemme DEFINITION 4.1. An {\bf   inextensional displacement} is defined as in the case of a single curved rod. It is an element of
the space $H^1_{\Gamma_0}({\cal S};\R^3)$ such that derivative tangential components vanish.

\noindent  We put
$$D_{I}=\Bigl\{U\in H^1_{\Gamma_0}({\cal S};\R^3)\ \ | \ \  {d U_i\over ds_i}\cdot\overrightarrow{T}_i =0,\quad
i\in\{1,\ldots,N\}\Bigr\}\leqno(15)$$  
\Lemme DEFINITION 4.2. An {\bf   extensional displacement} is a displacement belonging 
 to the orthogonal of $D_{I}$ in $H^1_{\Gamma_0}({\cal S};\R^3)$. 

\noindent The set of all  extensional displacements is denoted $D_{E}$.  We equip this  space  with the semi-norm 
$$||U||_{E}=\sqrt{\displaystyle\sum_{i=1}^N\int_0^{L_i}\Bigl|{dU_i\over ds_i}\cdot 
\overrightarrow{T}_i\Bigr|^2}\leqno(16)$$
\Lemme LEMMA 4.3.  The semi-norm $||.||_{E}$ is a norm equivalent to the $H^1_{\Gamma_0}({\cal S};\R^3)$-norm.

\noindent{\bf Proof : } It is clear that $||.||_{E}$ is a norm. For each $ i\in\{1,\ldots,N\}$, the range of the operator
$$U\in D_{E}\longmapsto \Bigl({dU_i\over ds_i}\cdot \overrightarrow{N}_i\, , \, {dU_i\over ds_i}\cdot
\overrightarrow{B}_i\Bigr)\in L^2(0,L_i;\R^2)$$ is a finite sub-space of $L^2(0,L_i;\R^2)$. Then we  prove by contradiction that
$||.||_{E}$ is a norm equivalent to the $H^1_{\Gamma_0}({\cal S};\R^3)$-norm.
\fin
 Let it be given a displacement $u$ belonging to $H^1_{\Gamma_0}({\cal S}_\delta;\R^3)$. The first component ${\cal U}$ of
the rods-structure elementary displacement
$U_e$ associated to $u$ is decomposed into a unique sum of an inextensional displacement and an extensional one.
$${\cal U}=U_I+U_E\qquad U_I\in D_{I},\quad U_E\in D_{E}\leqno(17)$$ 
\noindent From Lemma 3.4 we get the following estimates:
$$||U_E||^2_{H^1({\cal S};\R^3)}+\delta^2||U_I||^2_{H^1({\cal S};\R^3)}+\displaystyle\sum_{i=1}^N
\Bigl\|{dU_{Ii}\over ds_i}-{\cal R}_i\land\overrightarrow{T}_i \Bigr\|^2_{L^2(0,L_i; \R^3)}\le {C\over\delta^2}{\cal E}(u,{\cal
S}_\delta)\leqno(18)$$
\noindent{\bf 5. Properties of the limit of a sequence of displacements}
\vskip 1mm
 We recall the definition of  the unfolding operator (see [5]).  Let  $w$ be  in $L^1(\omega_{\delta, i})$. We denote 
${\cal T}_\delta( w)$ the function  belonging to $L^1(\omega_i)$ defined  as follows:
$${\cal T}_\delta (w )(s_i,Y_2,Y_3)=w(s_i,\delta Y_2,\delta Y_3),\qquad
\forall(s_i,Y_2,Y_3) \in \omega_i.  $$

 Let $u^\delta$ be a sequence of displacements of $H^1_{\Gamma_0}({\cal S}_\delta;\R^3)$ verifying
$${\cal E}(u^\delta,{\cal S}_\delta)\le C\delta^2\leqno(19)$$ 

 Displacement $u^\delta$ is decomposed into the sum of an elementary rods-structure  displacement and a residual one. The first
component of the elementary displacement is decomposed into the sum of an inextensional displacement and an extensional one.
From $(19)$ and estimates $(18)$, we deduce that it is possible to extract from these various sequences some subsequences,
still denoted in the same way, such that (see Proposition 7.1 of  [5])
$$\left\{\eqalign{
\delta {\cal U}^\delta,\ \ \delta U^\delta_I&\rightharpoonup U_I\quad\hbox{weakly in}\quad H^1({\cal S};\R^3)\cr 
\delta{\cal R}^\delta&\rightharpoonup {\cal R}\quad\hbox{weakly in}\quad H^1({\cal S};\R^3)\cr U^\delta_E&\rightharpoonup
U_E\quad\hbox{weakly in}\quad H^1({\cal S};\R^3)\cr
\delta{\cal T}_\delta(u^\delta_{|{\cal P}_{\delta,i}}) &\rightharpoonup U_{Ii}\quad\hbox{weakly in}\quad
H^1(\omega_i;\R^3)\cr {\cal T}_\delta(u^\delta_{|{\cal P}_{\delta, i}})-U^\delta_{Ii} &\rightharpoonup U_{Ei}-\Bigl(Y_2{dU_{Ii}
\over ds_i}\cdot \overrightarrow{N}_i+Y_3{dU_{Ii}\over ds_i}\cdot\overrightarrow{B}_i\Bigr)
\overrightarrow{T}_i\cr &\hskip 0.5cm
+\Theta_{U_I,i}\Bigl(-Y_3\overrightarrow{N}_i+Y_2\overrightarrow{B}_i\Bigr)\quad\hbox{weakly in}\quad
H^1(\omega_i;\R^3)\cr}\right.\leqno(20)$$ where  $\Theta_{U_I,i}={\cal R}_i\cdot\overrightarrow{T}_i$ is the limit rod torsion
angle of  ${\cal P}_{\delta,i}$. From estimates $(18)$ we derive also
$$\displaystyle{dU_{Ii}\over ds_i}={\cal R}_i\land \overrightarrow{T}_i,\qquad i\in\{1,\ldots,N\}\leqno(21)$$
\noindent After passing to the  limit there is a coupling between the inextensional displacement $U_I$ and the rod torsion angles
$\Theta_{U_I,i}$. They shall be reunited in a single functional space since we  have in any knot
$A$
$${dU_{Ii}\over ds_i}(a_i)={\cal R}(A)\land\overrightarrow{T}_i(a_i)\qquad i\in\{1,\ldots,N\}\leqno(22)$$ The condition 
obtained in $A$ expresses the rigidity of the junction in the neighborhood of $A$. The angles between two arcs joining in $A$
remain constant.
\medskip
  We consider now a new and last space containing the couple $(U_I,{\cal R})$. We set
$${\cal D}_I=\Bigl\{(V,{\cal A})\in H^1_{\Gamma_0}({\cal S};\R^3)
\times  H^1_{\Gamma_0}({\cal S};\R^3)\; | \; {dV_i\over ds_i}= {\cal A}_i\land\overrightarrow{T}_i,\; i\in\{1,\ldots,N\}
\Bigr\}\leqno(23)$$  
\noindent  where ${\cal D}_I$ is equipped with the norm
$$||(V,{\cal A})||=\sqrt{\sum_{i=1}^N\Bigl\|{d{\cal A}_i\over ds_i}\Bigr\|^2_{L^2(0,L_i;\R^3)}}\leqno(24)$$ This norm turns
${\cal D}_I$ into a Hilbert space. We denote $\Theta_{V,i}={\cal A}_i\cdot\overrightarrow{T}_i$ the rod torsion angles
associated to the couple $(V,{\cal A})$.
\vskip 1mm
\noindent{\ggras 6. Asymptotic behavior of structures made of curved rods}  
\vskip 1mm
\noindent The curved rods material is homogeneous and isotropic.  Let in ${\cal S}_\delta$ be the elasticity system
$$\left\{\eqalign{ - {\partial \over \partial x_j}\{\ a_{ljkh}{\partial u^\delta_k \over \partial x_h}\} & =F^\delta_l \enskip\hbox {
in }\enskip{\cal S}_\delta, \cr 
 u^\delta & =0 \hskip 1.4em \hbox { on }\hskip 0.7em\Gamma_{0\delta},\cr a_{ljkh} {\partial u^\delta_k \over\partial x_h}n_j &
=0\enskip\hbox { on }\hskip 0.7em  \Gamma_\delta,\qquad\Gamma_\delta=\partial {\cal S}_\delta \setminus
\Gamma_{0\delta},\cr}\right.\leqno(25)$$  where  the elasticity coefficients $a_{ljkh}$ are defined
by $$a_{ljkl}=\lambda\delta_{lj}\delta_{kh}+ \mu(\delta_{lk}\delta_{jh}+\delta_{lh}\delta_{jk}).$$  The constants
$\lambda$ and $\mu$ are the material Lam\'e coefficients. Obviously,  the coefficients  $a_{ljkh}$ satisfy  the hypothesis  of
coerciveness, i.e., there exists a constant $C_0>0$ such that
$$a_{ljkh}~\beta_{lj}~\beta_{kh}\geq C_0~\beta_{lj}~\beta_{lj},\quad\hbox{ for any }\;
\beta=(\beta_{lj})^{\vphantom{*}}_{{l,j}}~\hbox{ with }\beta_{lj}=\beta_{jl}.$$
 The applied forces $F_\delta$ are the sum of two different kinds of forces, 

$\bullet$ forces applied on each rod
$$\delta F_I+F_E,\qquad F_I, \, F_E\in L^2({\cal S} ; \R^3)\leqno(26)$$ 

$\bullet$ forces applied to the knots
$$ \bigl(\delta f_{A,I}+f_{A,E}\bigr){3\over 4\delta}\hbox{\bf 1}_{B(A;\delta)} \qquad \forall A\in {\cal N}\leqno(27)$$
\noindent This force $(27)$ is constant in the ball centered in $A$ and of radius $\delta$.  As in the case of a single rod, forces
$F_E$ and $f_{A,E}$ must not be concerned by the inextensional displacements in order to be able to work on the extensional
ones. Hence these forces satisfy the orthogonality condition
$$\displaystyle\int_{\cal S}F_E\cdot V+\sum_{A\in{\cal N}}f_{A,E}\cdot V(A)=0\quad\hbox{ for any
 }\ \ V\in D_I.\leqno(28)$$ According to estimates $(13)$ and $(18)$ we obtain 
$$\left\{\eqalign{ &\Bigl|{1\over\pi\delta^2} \int_{{\cal S}_\delta}F_\delta\cdot u-\int_{\cal S}F_E\cdot U_E -\int_{\cal
S}\delta F_I\cdot U_I +\sum_{A\in{\cal N}}\delta f_A\cdot U_I(A)\cr &+\sum_{A\in{\cal N}}f_{A,E}\cdot U_E(A)\Bigr|\le
C\sqrt{{\cal E}(u,{\cal S}_\delta)},\qquad \forall u\in H^1_{\Gamma_0}({\cal S}_\delta; \R^3).\cr}\right.\leqno(29)$$ 
\noindent The variational formulation of  problem $(25)$ is
$$\left\{\eqalign{
 u^\delta & \in H^1_{\Gamma_0}({\cal S}_\delta;\R^3)\cr 
\int_{{\cal P}_\delta} &\sigma_{kh}(u^\delta)\gamma_{kh}(v)=\int_{{\cal P}_\delta}F^\delta\cdot v\qquad \quad 
 \forall v  \in H^1_{\Gamma_0}({\cal S}_\delta;\R^3)\cr}\right.\leqno(30)$$
\noindent where $\sigma_{kh}(u^\delta)=a_{ljkh}\gamma_{lj}(u^\delta)$ are the stress tensor components.

 We can now estimate the solution $u^\delta$ to problem $(30)$. Due to $(29)$ we obtain
$\displaystyle{\cal E}(u^\delta,{\cal S}_\delta)\le C\delta^2$. Hence we can pass to the limit within the linearized elasticity
problem $(30)$.    
\vskip 1mm
\Lemme THEOREM 6.2. Extensional displacement $U_E$ is solution to the variational problem
$$\left\{\eqalign{ &U_E\in D_{E}\cr E&\sum_{i=1}^N\int_0^{L_i}{dU_{Ei}\over ds_i}\cdot\overrightarrow{T}_i{dV_i\over
ds_i}\cdot\overrightarrow{T}_i =\int_{\cal S}F_E\cdot V+\sum_{A\in{\cal N}}f_{A,E}\cdot V(A)
\qquad\quad\forall V\in D_{E}\cr }\right.\leqno(31)$$ 
\noindent The couple $(U_I,{\cal R})$ is solution to the variational problem
$$\left\{\eqalign{ 
&(U_I,{\cal R})\in {\cal D}_I,\cr   
& {E\over 3}\sum_{i=1}^N\int_0^{L_i}\hskip-1mm\Bigl\{\hskip-0.7mm{d^2U_{Ii}\over
ds^2_i}\cdot\overrightarrow{N}_i{d^2V_i\over ds^2_i}\cdot\overrightarrow{N}_i+\Bigl[{d^2U_{Ii}\over ds^2_i}
\cdot\overrightarrow{B}_i -c_i\Theta_{U_I,i}\Bigr]\hskip-1mm\Bigl[{d^2V_i\over
ds^2_i}\cdot\overrightarrow{B}_i-c_i\Theta_{V,i}\Bigr]\hskip-0.5mm\Bigl\}\cr    
+ &{\mu\over 3}\sum_{i=1}^N\int_0^{L_i}\hskip-0.7mm\Bigl[{d\Theta_{U_I,i}\over ds_i}+ c_i{dU_{Ii}\over
ds_i}\cdot\overrightarrow{B}_i\Bigr]\hskip-0.7mm\Bigr[{d\Theta_{V,i}\over ds_i}+c_i{dV_i\over
ds_i}\cdot\overrightarrow{B}_i\Bigr]\cr 
 =& \int_{\cal S}F_I\cdot V+\sum_{A\in{\cal N}}f_A\cdot V(A)\qquad\qquad
\forall(V,{\cal A})\in{\cal D}_I.\cr}\right.\leqno(32)$$

\noindent{\ggras Proof : }  The convergence results obtained in [5] can easily be transposed here since the elementary
rods-structure  displacement given in Lemma 3.4 and the elementary rod displacements defined by $(6)$ are very close according
to estimate $(12)$. Hence we obtain  the following limits of the stress-tensor components of sequence
$u^\delta$ in each reference domain $\omega_i$ (see  [5]):
$${\cal T}_\delta(\sigma_{kh}(u^\delta)_{|_{{\cal P}_{\delta,i}}})\rightharpoonup
\sigma_{i,kh}\quad\hbox{weakly in } L^2(\omega_i),$$ where 
$$\eqalign{  &\sigma_{i,11}=E\Bigl[{dU_{E,i}\over ds_i}\cdot\overrightarrow{T}_i-Y_2\Bigl({d^2U_{I,i}\over
ds^2_i}\cdot\overrightarrow{N}_i\Bigr) -Y_3\Bigl({d^2U_{I,i}\over
ds^2_i}\cdot\overrightarrow{B}_i-c_i\,\Theta_{U_I,i}\Bigr)\Bigr],\cr  &\sigma_{i,12}=\sigma_{i,21} =-{\mu Y_3\over
2}\Bigl[c_i\Bigl( {dU_{I,i}\over ds_i}\cdot\overrightarrow{B}_i\Bigl)+{d\Theta_{U_I,i}\over ds_i}\Bigr],\cr
&\sigma_{i,13}=\sigma_{i,31} =  \;\; {\mu Y_2\over 2}\Bigl[c_i\Bigl( {dU_{I,i}\over
ds_i}\cdot\overrightarrow{B}_i\Bigl)+{d\Theta_{U_I,i}\over ds_i}\Bigr],\cr 
&\sigma_{i,22}=\sigma_{i,33}=\vphantom{\Bigl[}\sigma_{i,23}=
\sigma_{i,32}=0,\cr}$$  and where $E=\displaystyle {\mu(3\lambda +2\mu)\over\lambda  + \mu}$  is the Young modulus.

 Let $V$ be an extensional displacement. We construct an admissible displacement $V^\delta$ of the whole structure ${\cal
S}_\delta$, by modifying $V$at the neighborhood of the junctions. Let us consider a knot $A$. This knot belongs to curve
$\gamma_i$ and its arc length is $a_i$. We modify $V$ in the neighborhood of 
$A$ by posing
$$V^\delta_i(s_i)=V_i(s_i)m\bigl({s_i-a_i\over\delta}\bigr)+V(A)
\Bigl(1-m\bigl({s_i-a_i \over\delta}\bigr)\Bigr)$$    Function $m$ is defined in $(10)$.
\noindent Displacement $V^\delta$  is an admissible displacement of the whole structure. Moreover
$$V^\delta\longrightarrow V\qquad \hbox{strongly in}\quad H^1({\cal S};\R^3).$$  
\noindent We take the test-displacement $V^\delta$ in $(30)$. We pass  to the limit and  we obtain  $(31)$.
\vskip 1mm  
\noindent Let us take now a couple $(V,{\cal A})\in {\cal D}_I$ and let us show $(32)$. We shall first construct a rods-structure
elementary displacement $v^\delta$  whose both components strongly converge in  $H^1_{\Gamma_0}({\cal S};\R^3)\times
H^1_{\Gamma_0}({\cal S};\R^3)$ towards $(V,{\cal A})$. Displacement $v^\delta$ is defined outside  the junctions by
$$v^\delta_i(s_i,y_2,y_3)={1\over\delta}V_i(s_i)+{\cal A}_i(s_i)\land\bigl({y_2\over\delta}
\overrightarrow{N}_i(s_i)+{y_3\over\delta}\overrightarrow{B}_i(s_i)\bigr\}$$   
\noindent Let $A$  be a knot belonging to curve $\gamma_i$ whose arc length is $a_i$. In the curve rod portion
$\Phi_i\bigl(]a_i-\rho\delta,a_i+\rho\delta[\times D(O;\delta)\bigr)$, the test-displacement $v^\delta$ is chosen rigid and is
given by
$$v^\delta(M)={1\over\delta}V(A)+{1\over\delta}{\cal A}(A)\land\overrightarrow{AM}\qquad M\in
\Phi_i\bigl(]a_i-\rho\delta,a_i+\rho\delta[\times D(O;\delta)\bigr)$$   
\noindent In the domains
$\Phi\bigl([a_i-2\rho\delta,a_i-\rho\delta]\times D(O;\delta)\bigr)$ and
$\Phi\bigl([a_i+\rho\delta,a_i+2\rho\delta]\times D(O;\delta)\bigr)$ displacement $v^\delta$ is defined as the rod elementary
displacement of Lemma 3.4 has been defined. Hence we have constructed an admissible displacement of the structure ${\cal
S}_\delta$. The unfolded strain tensor components  
 of sequence $v^\delta$ converge strongly in $ L^2(\omega_i)$ (see the proof of Theorems 7.1 and 7.2 of [5]).  Lastly we
consider $v^\delta$  as test-displacement in problem $(30)$. After having passed to the limit we obtain
$(32)$ with the couple $(V,{\cal A})$.\fin    
\noindent{\bf Corollary : } Variational problems $(31)$ and $(32)$ are coercive. Hence the whole sequences  converge toward
their limit. Following the same pattern as in Remark 7.4 of [5]  we prove that the convergences $(20)$ are strong. \fin  
\vskip 1mm
\noindent{\ggras 7. Complements}  
\vskip 1mm 
\noindent The extensional displacement $U_E$, solution to problem $(31)$, is $H^2$ class between two knots of a single arc.
Hence, at each knot $A$ we have
$$E\sum_{i=1}^N\bigl[\bigl({dU_{Ei}\over ds_i}(a_i+)-{dU_{Ei}\over ds_i}(a_i-)\bigr)\cdot
\overrightarrow{T}_i(a_i)\bigr]\overrightarrow{T}_i(a_i)=f_{A,E}\leqno(33)$$ This equality is the knots equilibrium law for a
structure made of curved rods. 
\medskip
\noindent{\bf REFERENCES  }
\medskip
\noindent [1] P.G. Ciarlet, H. Le Dret and R. Nzengwa. Junctions between three dimensional and two dimensional linearly elastic
structures.  J. Math. Pures Appl. 68 (1989), 261--295.

\noindent [2] D. Cioranescu, A. Damlamian and  G. Griso, Periodic unfolding and homogenization, C. R. Acad. Sci. Paris, Ser. I 335
(2002), 99--105.

\noindent [3] G. Griso, Etude asymptotiques de structures r\'eticul\'ees minces. Th\`ese Universit\'e Pierre et Marie Curie (Paris
VI), 1995.

\noindent [4] G. Griso, Asymptotic behavior of structures made of plates.  C. R. Acad. Sci. Paris, Ser. I 336 (2003), 101--107.

\noindent [5] G. Griso, Asymptotic behavior of curved rods by the unfolding method. (To appear in Math. Models Methods Appl.
Sci. )

\noindent [6] R. Jamal and E. Sanchez-Palencia, Th\'eorie asymptotique des tiges courbes anisotropes, C. R. Acad. Sci. Paris, Ser. I
332 (1996),  1099--1107.

\noindent [7] M. Jurak and J. Tamba\v ca, Derivation and justification of a curved rod model. Math. Models Methods Appl. Sci. 9, 7
(1999), 991--1015.

\noindent [8] H. Le Dret, Modeling of the junction between two rods, J. Math. Pures Appl. 68 (1989), 365--397.

\noindent [9] H. Le Dret. Modeling of a folded plate, Comput.  Mech., 5 (1990), 401--416.  

\noindent [10] J. Sanchez-Hubert and E. Sanchez-Palencia, Statics of curved rods on account of torsion and flexion. European Jour.
Mechanics A/Solids, 18 (1999), 365--390.

\noindent [11] E. Sanchez-Palencia, Sur certains probl\`emes de couplage de plaques et de barres. Equations aux d\'eriv\'ees
partielles et applications. Articles d\'edi\'es \`a J.L. Lions, Gauthier-Villars (1998), 725--745.
\bye